\documentclass[a4paper,11pt]{article}

\usepackage{amsmath,latexsym,amssymb}

\newtheorem{thm}{Theorem}
\newtheorem{prop}[thm]{Proposition}
\newtheorem{lemma}{Lemma}
\newtheorem{cor}{Corollary}

\newcommand{\proof}[1][]{{\it Proof#1: }}
\newcommand{\qed}[1][5mm]{\hspace*{\fill} $\Box$ \vspace{#1}}

\newcommand{\CC}{{\mathbf C}}

\newcommand{\del}{\partial}
\newcommand{\ON}{\operatorname}

\renewcommand{\a}{\alpha}
\renewcommand{\b}{\beta}
\newcommand{\cg}{\gamma}

\newcommand{\e}{\varepsilon}

\renewcommand{\l}{\lambda}

\newcommand{\s}{\sigma}

\newcommand{\inv}{^{^{-1}}}

\newcommand{\afami}{{\cal A}}
\newcommand{\cfami}{{\cal C}}

\newcommand{\sfami}{{\cal S}}

\newcommand{\labell}[1]{\label{#1}}

\newcommand{\br}{\ON{Br}}

\newcommand{\un}{\underline}

\begin{document}

\begin{center}
\Large\sc
Hurwitz stabilisers of some short redundant Artin systems for the braid
group $\br_3$\\[1.1cm]
\end{center}
\normalsize
{\bf \noindent
Michael\ L\"onne},\rm Institut f\"ur Mathematik, Am
Welfengarten 1, 30167 Hannover, Germany (e-mail: loenne@math.uni-hannover.de)

\begin{abstract}
We investigate the Hurwitz action of the braid group $\br_n$ on the $n$-fold
Cartesian product $\br_3^n$ and determine some stabilisers of its Artin systems.

Our algebraic result is complemented by a geometric study of families of plane
polynomial coverings of degree $3$. Together they lead to characterisations of
the set of {\it paths realised by degenerations} of the polynomials as defined
by Donaldson \cite{d}.
\end{abstract}

\section*{introduction}

The braid group $\br_n$ on $n$ strands has first been considered by A.\ Hurwitz
\cite{hu} and E.\ Artin \cite{ar}. Since then it has lead a prolific life in
such
diverse areas of mathematics as topology, combinatorics and algebraic
geometry.\\
In his original approach A. Hurwitz investigates the action of $\br_n$ on the
$n$-fold Cartesian product of a symmetric group. This action can formally be
defined on $G^n$ for any group $G$ and is determined by the action of the
standard generators
$$
\s_i(g_1,...,g_n)=(g_1,...,g_{i-1},g_ig_{i+1}g_i\inv,g_i,g_{i+2},...,g_n)
$$

Braid monodromy factorisations up to Hurwitz equivalence have proved
to be at the heart of the study of complex surfaces and symplectic
four-manifolds, cf.\ \cite{mo,mt,ak,kt}.

On the other hand people looked at stabilisers of Coxeter systems: Given a
symmetric matrix $M$ of positive integers $m_{ij}$ there is the associated
Coxeter group
$$
C=C_M:=\langle s_1,...,s_n\,|\,s_i^2=(s_is_j)^{m_{ij}}=1\rangle.
$$
The Hurwitz stabiliser $S_\cfami$ of the distinguished element
$\cfami=(s_1,...,s_n)\in C^n$ is shown to be 
$$
E_M:=\left\langle
e_{ij}:=\s_{j-1}\cdots\s_{i+1}\s_i^{m_{ij}}\s_{i+1}\inv\cdots\s_{j-1}\inv,
1\leq i<j\leq n\right\rangle
$$
in several important cases, \cite{bw,klui,cw,do}.

And there is always a geometric flavour to these results; Birman and Wajnryb
\cite{bw} applied their results to get a presentation for the mapping class
group of surfaces with at most one boundary component and the stabilisers of the
classical Coxeter systems for the finite Coxeter groups
investigated by Catanese/Wajnryb (type $A$) and D\"orner (type $ADE$) had been
shown to be the fundamental groups of the bifurcation complements of the
corresponding simple singularities, cf.\ \cite{lo}.
\\

On the geometric side our research was initiated by the article \cite{d} in
'Mathematics: frontiers and perspectives'.

Given a plane polynomial $g_1=g_1(x,y)$ the projection to the $x$-axis defines
a function $f_1$ on the plane curve $C_1$ defined by $g_1$.
Suppose now $C_1$ is smooth and $f_1$ simply branched, i.e.\ $f_1$ is a Morse
function with singular values $b_0,...,b_n\in\CC$, which are called the branch
points.
Then Donaldson asks for the set of paths which can be realised by deformations
of $g_1$, cf. \cite[ p.56]{d}.
We call such paths {\it vanishing arcs}, since we have the following
characterisation:

An embedded path $\cg$ with endpoints $b_0,b_1$ and disjoint from the branch
set otherwise is a vanishing arc, if there is a smooth family $g_t,t\in[0,1]$
of polynomials and a homotopy $H=H(s,t)$ such that 
\begin{enumerate}
\item
the curves $C_t$ are smooth and $f_t$ simply branched for $t>0$,
\item
the curve $C_0$ is smooth except for a single ordinary double point $P$ and $f_0$
is simply branched on the smooth locus and injective on its critical points,
\item
$H$ is a homotopy from $\cg$ to the constant map to $f_0(P)$ such that $H_t$
meets the branch set of $f_t$ in $H_t(\{0,1\})$.
\end{enumerate}
So a vanishing arc can be contracted to point relative to the branch locus.
\\ 

We first review results on Coxeter systems and tailor them to our needs. In
subsequent sections we investigate stabilisers of Artin systems
$(a_1,...,a_n)\in A^n$ which are Artin analogues of Coxeter systems and
associated to Artin groups
$$
A = A_{M}:=\langle
a_1,...,a_n\,|\,\underbrace{a_ia_ja_i...}_{m_{ij}\text{
factors}}=
\underbrace{a_ja_ia_j...}_{m_{ij}\text{ factors}}\rangle.
$$
To handle degenerations of plane polynomials we choose the set up of polynomial
covers introduced by Hansen \cite{ha} and define the notion of bifurcation
braid monodromy.
The paper closes with a characterisation of sets of vanishing arcs in the
spirit of \cite{d}.
\\

Please note that while we do the computations for all the Hurwitz actions in
the proofs, we only state a lot of not so obvious braid relations which should
be checked by hand or using a symbolic computation package in case of doubt.

We use $\s^\b$ and $(\s)\b$ to denote the action of a braid $\b$ on a braid
$\s$ by conjugation on the right.

\section*{stabilisers of Coxeter systems}

We start with the Hurwitz action on Cartesian powers of the symmetric group
$\sfami_3$ on three elements. The Coxeter presentation for $\sfami_3$ is
generated by two transpositions $s$ and $t$:
$$
\sfami_3\,\cong\,C:=\langle s,t\,|\,s^2=t^2=(st)^3=1\rangle.
$$
Any non-constant $n$-tuple $(s_1,...,s_n)$ with $s_i\in\{s,t\}$ is not only an
element in $C^n$ but a Coxeter system with Coxeter group $C$ and matrix
$$
M:=\big(m_{ij}\big)_{1\leq i,j\leq n},\text{ with } m_{ij}=\left\{
\begin{array}{lll}
3 & \text{if} & s_i\neq s_j\\
1 & \text{if} & s_i=s_j
\end{array}\right.
$$

For some of these Coxeter systems the stabiliser group can be found with
results
of Birman and Wajnryb.

\begin{lemma}
\labell{bw-gen}
The stabiliser group of the Coxeter systems $(s,t,t,...,t)$ and
$(s,s,t,...,t)$ can
be given by
\begin{eqnarray*}
S_C & = & \langle
\s_1^3,\s_2,...,\s_{n-1},
(\s_4)\s_3\s_2\s_1^2\s_2\s_3^2\s_2\s_1 \rangle,\\[4mm]
\text{respectively}\\[3mm]
S_C & = & \langle
\s_1,\s_2^3,\s_3,...,\s_{n-1},
(\s_5)\s_4\s_3\s_2^2\s_3\s_4^2\s_3\s_2,
(\s_3)\s_2\inv\s_1^{-2}\s_2^2\s_1 \rangle.\\
\end{eqnarray*}
\end{lemma}

\proof
The first claim is actually proved in \cite{bw}. The argument can then be
extended to cover the other case, cf.\ \cite{klui}.
\qed

Starting from this lemma we can get the stabiliser groups of the alternating
Coxeter systems $\cfami_n:=(s,t,s,t,...)\in C^n$ with matrix
$$
\label{matrix}
M_n:=\big(m_{ij}\big)_{1\leq i,j\leq n},\text{ with } m_{ij}=\left\{
\begin{array}{lll}
3 & \text{if} & i\not\equiv j \mod 2\\
1 & \text{if} & i\equiv j \mod 2
\end{array}\right.
$$

\begin{lemma}
\labell{conj-elm}
The stabiliser groups $S_{\cfami_n}$ of the alternating Coxeter systems
of length $n$ are conjugated to the stabiliser group of a
Coxeter system $(s,t,t,...,t)$ or $(s,s,t,...,t)$ of equal length. For
small $n$
we have:
$$
\begin{array}{llrcl}
S_{(s,t,s)} & = & \s_1\inv &S_{(s,t,t)}& \s_1\\
S_{(s,t,s,t)} & = & \s_2\inv\s_1 &S_{(s,t,t,t)}& \s_1\inv\s_2\\
S_{(s,t,s,t,s)} & = & \s_3\inv\s_2\s_1\inv &S_{(s,t,t,t,t)}& \s_1\s_2\inv\s_3\\
S_{(s,t,s,t,s,t)} & = & \s_4\inv\s_3\s_2\inv &S_{(s,s,t,t,t,t)}&
\s_2\s_3\inv\s_4\\
\end{array}
$$
\end{lemma}

\proof
The general claim is immediate from the observation that the alternating
system is in the same $\br_n$ orbit as at least one of the reference systems.\\
For the second claim we compute the action of the conjugating element on the
alternating system of two transpositions $s,t$ in the permutation group of $3$
elements, with $r$ denoting the third transposition.
$$
\begin{array}{rcr@{}l}
\s_1  (\un{s,t},s) &
= & & (r,s,s)\\[3mm]
\s_1\inv\s_2  (s,\un{t,s},t) &
= & \s_1\inv & (\un{s,r},t,t)\\
& = & & (r,t,t,t)\\[3mm]
\s_1\s_2\inv\s_3  (s,t,\un{s,t},s) &
= & \s_1\s_2\inv & (s,\un{t,r},s,s)\\
& = & \s_1 & (\un{s,r},s,s,s)\\
& = & & (t,s,s,s,s)\\[3mm]
\s_2\s_3\inv\s_4  (s,t,s,\un{t,s},t) &
= & \s_2\s_3\inv & (s,t,\un{s,r},t,t)\\
& = & \s_2 & (s,\un{t,r},t,t,t)\\
& = & & (s,s,t,t,t,t)\\[3mm]
\end{array}
$$
Surely the stabiliser does not depend on the choice of two non-commuting
transpositions and hence we are done.
\qed

Next we show that the $E_n:=E_{M_n}$ is a set of stabilisers and that its set
of generators can easily be reduced.

\begin{lemma}
\labell{E-gen}
The group $E_n$ is generated by $\s_1^3,\s_{i+1}\s_i\s_{i+1}\inv$,
$i=1,...,n-2$ and is contained in the stabiliser subgroup $S_{\cfami_n}$ of the
alternating Coxeter system $\cfami_n$.
\end{lemma}

\proof
First note that the given generators are just $e_{12}, e_{i,i+2}$,
$i=1,...,n-2$. For the remaining $e_{ij}$ we have the following relations:
$$
\begin{array}{rcll}
e_{ij}&=&(e_{i,i+2})e\inv_{i+2,i+4}e\inv_{i+4,i+6}...e\inv_{j-2,j}, &
i\equiv j\mod 2,\\
e_{ij}&=&(e_{12})e_{13}e_{24}...
e_{i-1,i+1}e\inv_{i+1,i+3}e\inv_{i+3,i+5}...e\inv_{j-2,j}, & i\not\equiv
j\mod 2.
\end{array}
$$
For the second claim it thus suffices to check that
$$
\begin{array}{r@{}lcr@{}lcr@{}lcl}
\s_1^3 & (s,t) & = & \s_1^2 & (r,s) & = & \s_1 & (t,r) & = & (s,t),\\
\s_2\s_1\s_2\inv & (s,\un{t,s}) & = & \s_2\s_1 & (\un{s,s},r) & = & \s_2 &
(s,\un{s,r}) & = & (s,t,s).
\end{array}
$$
\pagebreak
\qed[0mm]

In fact we can choose two additional stabilisers to generate the stabiliser
group $S_{\cfami_n}$. Our choice
$\tau_1=(\s_1)\s_2\s_3\inv\s_4,\tau_2=(\s_2)\s_3\s_4\inv\s_5$ is motivated by
their action on Artin systems, cf.\ lemma \ref{tauArtin}.

\begin{lemma}
\labell{alt-gen}
The stabiliser groups $S_{\cfami_n}$ of the alternating Coxeter systems
of lengths $n=3,4,5,6$ can be given as follows:
$$
\begin{array}{lclllcl}
S_{(s,t,s)} & = & E_3,  & &
S_{(s,t,s,t)} & = & E_4,\\
S_{(s,t,s,t,s)} & = & \langle E_5,\tau_1\rangle,  & &
S_{(s,t,s,t,s,t)} & = & \langle
E_6,\tau_1,\tau_2\rangle.
\end{array}
$$
\end{lemma}

\proof
To show that $\tau_1,\tau_2$ belong to the stabilisers of $\cfami_5,\cfami_6$
resp.\, it obviously suffices to show $\tau_1\in S_{\cfami_5}$:
$$
\begin{array}{crl}
& \s_4\inv\s_3\s_2\inv\s_1\s_2\s_3\inv\s_4 & (s,t,s,\un{t,s})\\
= & \s_4\inv\s_3\s_2\inv\s_1\s_2\s_3\inv & (s,t,\un{s,r},t)\\
= & \s_4\inv\s_3\s_2\inv\s_1\s_2 & (s,\un{t,r},t,t)\\
= & \s_4\inv\s_3\s_2\inv\s_1 & (\un{s,s},t,t,t)\\
= & \s_4\inv\s_3\s_2\inv & (s,\un{s,t},t,t)\\
= & \s_4\inv\s_3 & (s,t,\un{r,t},t)\\
= & \s_4\inv & (s,t,s,\un{r,t})\\
= & & (s,t,s,t,s)
\end{array}
$$
So the given groups are shown to be stabilising.

To prove the inverse implication we note that the stabiliser groups are
generated by the following elements which are obtained from generators of lemma
\ref{bw-gen} using conjugation as provided by lemma \ref{conj-elm}:
$$
\begin{array}{ll}
n=3: &
 \s_1^3,(\s_2)\s_1 ,\\
n=4: &
 (\s_1^3)\s_2,(\s_2)\s_1\inv\s_2,(\s_3)\s_2 ,\\
n=5: &
 (\s_1^3)\s_2\inv\s_3,(\s_2)\s_1\s_2\inv\s_3,
(\s_3)\s_2\inv\s_3,(\s_4)\s_3,\\
& \hspace{7mm}(\s_4)\s_3\s_2\s_1^2\s_2\s_3^2\s_2\s_1\s_1\s_2\inv\s_3,\\
n=6: &
 (\s_1)\s_2\s_3\inv\s_4,(\s_2^3)\s_3\inv\s_4,(\s_3)\s_2\s_3\inv\s_4,
(\s_4)\s_3\inv\s_4,(\s_5)\s_4,\\
& \hspace{7mm}(\s_5)\s_4\s_3\s_2^2\s_3\s_4^2\s_3\s_2\s_2\s_3\inv \s_4,
(\s_3)\s_2\inv\s_1^{-2}\s_2^2\s_1\s_2\s_3\inv\s_4.
\end{array}
$$
The claim now follow from the fact that all these elements can be expressed by
elements of the given groups:
$$
\begin{array}{lrlrlrl}
n=3: && \s_1^3 &&(\s_2)\s_1\\
& \hspace{1mm}= & e_{12}& \hspace{6mm}= & e_{13}\\[1mm]
n=4: && (\s_1^3)\s_2 && (\s_2)\s_1\inv\s_2 && (\s_3)\s_2\\
& = & (e_{12})e_{23}\inv e_{13}\inv & = & (e_{13})e_{23}
& \hspace{6mm}= & e_{24}\\[2mm]
n=5: && (\s_1^3)\s_2\inv\s_3 && (\s_2)\s_1\s_2\inv\s_3 && (\s_3)\s_2\inv \s_3\\
& = & (e_{34})e_{13}\inv & = & (e_{24})e_{34}e_{13}\inv & = & (e_{24})e_{34}
\\[1mm]
&& (\s_4)\s_3 &&
\multicolumn{3}{l}{(\s_4)\s_3\s_2\s_1^2\s_2\s_3^2\s_2\s_1\s_1\s_2\inv\s_3}\\
& = & e_{35} & = & 
\multicolumn{3}{l}{(\tau_1)e_{35}e_{45}e_{34}e_{15}\inv}\\[2mm] 
n=6: && (\s_2^3)\s_3\inv\s_4 && (\s_3)\s_2\s_3\inv\s_4 && (\s_2)\s_3\inv\s_4\\ 
& = & (e_{45})e_{24}\inv & = & (e_{35})e_{45}e_{24}\inv & = & (e_{35})e_{45}
\\[1mm]
&& (\s_5)\s_4 &&
\multicolumn{3}{l}{(\s_5)\s_4\s_3\s_2^2\s_3\s_4^2\s_3\s_2\s_2\s_3\inv\s_4}\\
& = & e_{46} & = & 
\multicolumn{3}{l}{(\tau_2)e_{46}e_{56}e_{45}e_{26}\inv}\\[1mm] 
&& (\s_1)\s_2\s_3\inv\s_4 &&
\multicolumn{3}{l}{(\s_3)\s_2\inv\s_1^{-2}\s_2^2\s_1\s_2\s_3\inv\s_4}\\
& = & \tau_1 & = & 
\multicolumn{3}{l}{e_{13}}
\end{array}
$$
\pagebreak
\qed

While the claim of the preceding lemma holds true for $n>6$ the next result
does not and is therefore special to $n=5,6$.

\begin{lemma}
\labell{normal}
In case $n=5,6$ the group $E_n$ is normal in the stabiliser group
$S_{\cfami_n}$.
\end{lemma}

\proof
It suffices to
show that conjugation by $\tau_1$ resp.\ $\tau_1,\tau_2$ maps generators of
$E_n$ to $E_n$:
$$
\begin{array}{llcllcl}
n=5,6: &
(e_{12})\tau_1\inv & = & (e_{45})e_{24}\inv, &
(e_{13})\tau_1\inv & = & (e_{13})e_{12}\inv e_{24} e_{45}e_{24}\inv\\
& (e_{24})\tau_1\inv & = & e_{24} &
(e_{35})\tau_1\inv & = & (e_{35})e_{15}e_{12}\inv e_{15}\inv
e_{12}\inv e_{45}\\[2mm]
n=6: &
(e_{23})\tau_2\inv & = & (e_{56})e_{35}\inv, &
(e_{24})\tau_2\inv & = & (e_{24})e_{23}\inv e_{35} e_{56}e_{35}\inv\\
& (e_{35})\tau_2\inv & = & e_{35} &
(e_{46})\tau_2\inv & = & (e_{46})e_{26}e_{23}\inv e_{26}\inv
e_{23}\inv e_{56}\\
&(e_{12})\tau_2\inv & = & (e_{56})e_{35}\inv e_{13}\inv, &
(e_{56})\tau_1\inv & = & (e_{12})e_{15}\inv e_{12}\inv e_{45}e_{46}

\end{array}
$$
\qed

\section*{stabilisers for redundant Artin systems}

In analogy with the previous section we start with a standard presentation for
the braid group on three strands generated by two braids $a$ and $b$:
$$
A:=\langle a,b\,|\,aba=bab\rangle.
$$
As we may expect the alternating Artin systems $\afami_n:=(a,b,a,b,...)$ of
length $n\geq2$ have matrix $M_n$ and Artin group $A$.

\begin{lemma}
\labell{Hact}
Let $A$ act by elementwise conjugation on the $n$-fold
Cartesian product $A^n$, then the stabiliser subgroup of the alternating Artin
system $\afami_n$ has trivial intersection with the subgroup $H=\langle
a^2,b^2\rangle$ if $n\geq2$.
\end{lemma}

\proof
By the positive solution of the Tits conjecture \cite{cp}, $H$ is freely
generated by
$a^2$ and
$b^2$. On the other hand if $n\geq2$ any element in the stabiliser
must actually
belong to the center of $A\cong\br_3$. But the intersection of the center of
$A$ with $H$
must be contained in the center of $H$ which is trivial since $H$ is free.
\qed

\begin{lemma}
\labell{tauArtin}
The braids $\tau_1,\tau_2$ act on elements of the $H$-orbit of $(a,b,a,b,a,b)$
by overall conjugation with $b^{-2}$ resp.\ $a^{-2}$.
\end{lemma}

\proof
Since overall conjugation commutes with the braid action it suffices to
prove the claim
for the action of $\tau_1,\tau_2$ on $(a,b,a,b,a,b)$,
even one of these cases
suffices by symmetry:
$$
\begin{array}{crl}
& \s_4\inv\s_3\s_2\inv\s_1\s_2\s_3\inv\s_4 & (a ,b,a ,\un{b,a} ,b)\\[2mm]
= & \s_4\inv\s_3\s_2\inv\s_1\s_2\s_3\inv & (a ,b,\un{a ,ba
b\inv},b,b)\\[2mm]
= & \s_4\inv\s_3\s_2\inv\s_1\s_2 & (a ,b,ba b\inv,ba \inv
b\inv
a ba b\inv,b,b)\\[2mm]
= & \s_4\inv\s_3\s_2\inv\s_1\s_2 & (a \un{,b,ba
b\inv},b,b,b)\\[2mm]
= & \s_4\inv\s_3\s_2\inv\s_1 & (\un{a ,bba b\inv
b\inv},b,b,b,b)\\[2mm]
= & \s_4\inv\s_3\s_2\inv & (a bba b\inv b\inv a \inv,\un{a
,b},b,b,b)\\[2mm]
= & \s_4\inv\s_3 & (a bba b\inv b\inv a
\inv,b,\un{b\inv a b,b},b,b)\\[2mm]
= & \s_4\inv & (a bba b\inv b\inv a \inv,b,b\inv
a bbb\inv a \inv b,\un{b\inv
a b,b},b)\\[2mm]
= & & (a bba b\inv b\inv a \inv,b,b\inv a ba
\inv b,b,b\inv b\inv a bb,b)
\end{array}
$$
And the last line equals $(b^{-2}a b^2,b ,b^{-2}a
b^2,b,b^{-2}a
b^2,b)$ since
$$
a bba b\inv b\inv a \inv=a ba \inv ba
b\inv a \inv=b\inv a bbb\inv a \inv b=
b\inv a ba \inv b=b\inv b\inv a bb.
$$
{}\qed

\begin{cor}
\label{corr}
The braid $\tau_1$ acts trivially on alternating Coxeter systems but
non-trivially on
alternating Artin systems for $n\geq5$.
\end{cor}

\begin{prop}
\labell{semi-prod}
The groups $S_{\cfami_5},S_{\cfami_6}$ are semi direct products of their normal
subgroups $E_5$ resp.\ $E_6$ and a free subgroup freely
generated by $\tau_1$ resp.\ $\tau_1,\tau_2$.
\end{prop}

\proof
The subgroup generated by $\tau_1$ resp.\ $\tau_1,\tau_2$ acts freely on
the $H$-orbit
of $(a ,b,a ,b,a )$ resp.\ $(a ,b,a ,b,a ,b)$,
and therefore is a free subgroup. The
claim is then immediate from the normality of the group $E_n$ in $S_{\cfami_n}$
for $n=5,6$, cf.\ lemma \ref{normal}.
\qed

\begin{thm}
Let $S_{\afami_n},S_{\cfami_n}$ be the stabilisers of the Artin resp. Coxeter
system associated with the
$n\times n$ matrix
$$
M_n:=(m_{ij}),\,1\leq i,j\leq n,\text{ with } m_{ij}=\left\{
\begin{array}{lll}
3 & \text{if} & i\not\equiv j \mod 2\\
1 & \text{if} & i\equiv j \mod 2
\end{array}\right.
$$
and denote by $E_n$ the subgroup of the braid group $\br_n$ generated by
elements
\linebreak
$e_{ij}:=\s_{j-1}\cdots\s_{i+1}\s_i^{m_{ij}}\s_{i+1}\inv\cdots\s_{j-1}\inv,
1\leq i<j\leq n$,
then the following relations hold:
$$
\begin{array}{cccccll}
E_n & = & S_{\afami_n} & = & S_{\cfami_n} & \,\,\,\text{ if} & n=2,3,4,\\
E_n & = & S_{\afami_n} & < & S_{\cfami_n} & \,\,\,\text{ if} & n=5,6,\\
E_n & \leq & S_{\afami_n} & < & S_{\cfami_n} & \,\,\,\text{ if} & n\geq7.
\end{array}
$$
\end{thm}

\proof
By straightforward calculation the inclusion $E_n\leq S_{\afami_n}$ is shown,
whereas the inclusion $S_{\afami_n}\leq S_{\cfami_n}$ is obvious.
Equality $E_n=S_{\cfami_n}$ for $n=2,3,4$ is shown in lemma \ref{alt-gen}.
Strict inclusion $S_{\afami_n}< S_{\cfami_n}$ for $n\geq5$ follows from
corollary
\ref{corr}. In case $n=5,6$ finally each braid in $S_{\cfami_n}$ can be
written as a product
$\tau e$ with
$\tau\in\langle\tau_1,\tau_2\rangle$ and $e\in E_n$ by prop.\
\ref{semi-prod}. Hence
it is immediate by lemma \ref{Hact} that $E_n$ is the stabiliser group
$S_{\afami_n}$ in these cases.
\qed

\section*{conjugacy classes of simple braids}

The aim of this section is to exhibit the set of simple braids in $E_n$ as a
single $E_n$ conjugation class for $n\leq6$, where a braid is called simple if
it is isotopic to a half twist associated to a path connecting two punctures.

We extend the conjugation action and its exponential notation to sets, so the
conjugation orbit of $e_{13}$ in $E_n$ is denoted by
$e_{12}^{E_n}:=\{e_{12}\}^{E_n}$.

\begin{lemma}
\labell{Scases}
For $n\leq6$ there is a set of braids $T_{B/S}$ such that $\br_n=T_{B/S}\cdot
S_{\cfami_n}$ and for all $\cg\in{T}_{B/S}$:
\begin{eqnarray*}
e_{13}^\cg\in E_n & \Rightarrow & e_{13}^\cg\in e_{13}^{E_n}.
\end{eqnarray*}
\end{lemma}

\proof
We defined $S_{\cfami_n}$ to be a stabiliser for the Hurwitz action of $\br_n$ on
the finite set $\sfami_3^n$, hence the stabiliser is of finite index. Our
strategy is to construct $T$ as a Schreier left transversal for the
$S_{\cfami_n}$-cosets.

First note that the cosets are in bijection to the elements of the
$\br_n$-orbit of the Coxeter system. Moreover we may assume that the Schreier
transversal contains only positive braids since $\s_i\inv$ and $\s_i^2$ act
the same way on all elements of $\sfami_3^n$.

Such transversal $T_{B/S}$ for $E_6$ can be found, e.g.\ by a short symbolic
computation, containing $240$ elements. Since $E_n= S_\afami$ we extract a
list of all elements $\cg$ such that $e_{13}^\cg$ stabilises $\afami_6$.

For such elements we get in fact $18$ elements $e_{13}^\cg$, which can be shown
to be conjugated in $E_n$ to elements $e_{i,i+2}$ and hence to $e_{13}$:
$$
\begin{array}{lll}
(e_{13})\s_3\s_1=(e_{24})e_{13}\inv &
(e_{13})\s_1\s_2=(e_{13})e_{12}\inv \\
(e_{13})\s_3\s_2\s_2=(e_{13})e_{23}e_{24}\inv &
(e_{13})\s_2\s_2\s_3=(e_{13})e_{23}e_{24} \\
(e_{13})\s_3\s_3\s_1\s_2=(e_{24})e_{12}e_{34} &
(e_{13})\s_3\s_4\s_1\s_1=(e_{24})e_{34}e_{35}e_{12} \\
(e_{13})\s_3\s_3\s_2\s_1\s_1=(e_{24})e_{12}\inv e_{24}\inv &
(e_{13})\s_3\s_2\s_2\s_3\s_1=(e_{24})e_{34}e_{13}^{-2} \\
(e_{13})\s_2\s_2\s_3\s_3\s_1=(e_{24})e_{23} &
(e_{13})\s_2\s_2\s_3\s_4\s_1=(e_{24})e_{34}e_{35} \\
(e_{13})\s_3\s_3\s_4\s_4\s_2=(e_{12})e_{23}e_{35}\inv e_{45} &
(e_{13})\s_3\s_3\s_4\s_5\s_2=(e_{13})e_{23}e_{35}\inv e_{45}e_{46} \\
(e_{13})\s_3\s_4\s_4\s_5\s_2=(e_{13})e_{23}e_{35}\inv e_{45}e_{46}e_{24} &
(e_{13})\s_3\s_4\s_4\s_2\s_1\s_1=(e_{35})e_{45}e_{13}\inv \\
(e_{13})\s_3\s_4\s_5\s_2\s_1\s_1=(e_{35})e_{45}e_{13}\inv e_{46} &
(e_{13})\s_3\s_4\s_2\s_3\s_3\s_1=(e_{24})e_{34}e_{35}\inv \\
(e_{13})\s_3\s_4\s_5\s_1\s_1\s_2=(e_{13})e_{12}e_{35}\inv e_{45}e_{13}
\inv e_{12}\inv e_{46}\\
(e_{13})\s_3\s_4\s_4\s_5\s_5\s_2\s_3\s_1=(e_{46})e_{56}e_{24}e_{35}e_{56}
\end{array}
$$
\qed

\begin{lemma}
\labell{tcases}
Let $\tau\in S_{\cfami_5}$, resp.\ $S_{\cfami_6}$ be a freely reduced word in
$\tau_1$, resp.\
$\tau_1,\tau_2$, then
$$
e_{13}^\tau\in e_{13}^{E_n}.
$$
\end{lemma}

\proof
For the at most four words of
length one we have:
$$
\begin{array}{lcl}
(e_{13})\tau_1 & = & (e_{13})e\inv_{35}e_{45}e_{35}e_{45}e_{13}e_{12}\inv \\
(e_{13})\tau_1\inv & = & (e_{13})e_{12}\inv e_{23}e_{45}e_{24}\inv\\
(e_{13})\tau_2 & = & e_{13}\\
(e_{13})\tau_2\inv & = & e_{13}
\end{array}
$$
To conclude by induction on the word length we next consider a word $\tau$ which
is the concatenation of a word
$\tau'$ and a letter $\tau''$. By induction hypothesis
$e_{13}^{\tau'}=e_{13}^e$ for some $e\in E_n$ and by the normality of
$E_n$ in $\langle E_n,\tau_1\rangle$, resp.\ $\langle
E_n,\tau_1,\tau_2\rangle$ we can write
$e\tau''=\tau''e'$ for a suitable $e'\in E_n$. Hence
$$
e_{13}^{\tau}=e_{13}^{\tau'\tau''}=e_{13}^{e\tau''}
=e_{13}^{\tau''e'}\in e_{13}^{E_n}.
$$
\qed

\begin{prop}
\labell{conclass}
The intersection of the conjugacy class of half twists in
$\br_n,\,n\leq6$, with $E_n$ coincides with the conjugacy class of $e_{13}$
in $E_n$.
\end{prop}

\proof
One inclusion is obvious. So we pick any $\b\in\br_n$ such that $e_{13}^\b\in
E_n$, which may be factorised as $\b=\cg\tau e$ with $\cg\in T_{B/S}$,
$\tau\in\langle\tau_1,\tau_2\rangle$, $e\in E_n$.
Since $\tau e$ normalises $E_n$, we get $e_{13}^\cg\in E_n$ and
$e_{13}^\cg=e_{13}^{e'}$ for some $e'\in E_n$ by lemma \ref{Scases}.

By lemma \ref{normal} and lemma \ref{tcases} there are $e'',e'''\in E_n$ such
that:
$$
e_{13}^{\b} =e_{13}^{\cg\tau e} =e_{13}^{e'\tau e} 
=e_{13}^{\tau e''e} =e_{13}^{e'''e''e} \in e_{13}^{E_n}.
$$
\qed

\section*{bundles and monodromy}

We owe the following exposition of polynomial covers to \cite{cs} and
\cite{ha}.

On a connected topological manifold $X$ a simple Weierstrass polynomial of
degree~$d$ is
a map $f:X\times\CC\to\CC$ given by
$$
f(x,z):=z^d+\sum_{i=1}^d c_i(x)z^{d-i},
$$
with continuous coefficient maps $c_i:X\to\CC$, and with no multiple roots
for any
$x\in X$.
Given such a function $f$ the first coordinate projection map onto $X$ may be
restricted
to the subspace
$$
Y_f:=\{(x,z)\in X\times\CC\,|\,f(x,z)=0\}
$$
defining a $d$-fold cover $\pi_f$ onto $X$, the polynomial cover associated
to $f$, or
to the complement $X\times\CC\setminus Y_f$ defining a fibre bundle over
$X$ with fibre
diffeomorphic to a $d$-punctured disc, the punctured disc bundle associated
to $f$.

A finite unramified cover is called polynomial if it is equivalent to a
polynomial
cover for some simple Weierstrass function as above. Any cover $\pi:Y\to
X$ gives rise to a monodromy homomorphism from $\pi_1(X,x)$ to the symmetric
group $\sfami(\pi\inv(x))$,
which serves for a natural characterisation of polynomial covers: 

\begin{prop}[\cite{ha}]
\labell{alg-lift}
An unramified cover of degree $d$ is  polynomial if and only if its monodromy
homomorphism to the symmetric group $\sfami_d$ lifts along the natural
homomorphism $\br_d\to\sfami_d$.
\end{prop}

There is a natural way to get from a Coxeter system of
length~$n$ of a
symmetric group $\sfami_d$ to a finite cover: Given an $n$-punctured disc and
a geometric
basis for its fundamental group, which is a choice of free generators
$\cg_1,...,\cg_n$ such that
\begin{enumerate}
\item
the generator are represented by disjoint paths freely homotopic to positive
loops around single punctures,
\item
the product $\cg_1\cdots\cg_n$ is freely homotopic to the positive boundary of
the disc.
\end{enumerate}

A homomorphism to $\sfami_d$ is obtained by assigning to
these generators the elements of the Coxeter system. The preimage
of any
subgroup isomorphic to $\sfami_{d-1}$ determines subgroups of the fundamental
group in a unique conjugacy class and thus a well defined finite cover of the
punctured disc.\\

The corresponding result associates with an Artin system of length~$n$ for
the braid
group $\br_d$ a $d$-punctured disc bundle, once a geometric basis for the
fundamental
group of an $n$-punctured disc has been chosen. Here the basis elements are
mapped to
the generators of the Artin system, so a homomorphism to the braid group is
obtained.\\
Since the space of monic polynomials is an Eilenberg-MacLane space for the
braid group,
there is a smooth classifying map for this homomorphism. Pulling back the
tautological
simple Weierstrass polynomial we get a simple Weierstrass polynomial on the
$n$-punctured disc, and the associated punctured disc bundle is the one we
aim for.

\section*{families of polynomials on the plane}

We enter now the realm of complex geometry where there is an abundance of
covers and
bundles as defined in the previous section.

\begin{description}
\item[Example 1:]
Given a branched cover $Y\to X$ of a complex manifold $X$, the restriction
to the
branch complement is a finite topological cover. Its monodromy is also
called the
monodromy of the branched cover.

\item[Example 2:]
Given a plane curve $C\in\CC^2$ and a projection $p:\CC^2\to\CC$ such that
$p|_C$ is a
finite branched cover, then restricted to the preimage of the complement of
the branch
locus $p|_C$ is a polynomial cover and $p|_{\CC^2-C}$ is a punctured disc
bundle. The
corresponding structure homomorphism to the braid group is called the braid
monodromy.

\end{description}

The second example leads a straight way to the following generalised notion
of braid
monodromy:

\begin{description}

\item[Definition:]
Given a divisor $D\subset T\times\CC$ such that the map $pr|_D$ induced
from the first
projection $pr=pr_1$ onto $T$ is a finite cover, the restriction of $pr$ to the
intersection of the complement of $D$ and the preimage of the branch
complement is a
punctured disc bundle and its structure homomorphism is called the
(generalised) braid
monodromy of $D$.

\end{description}

In favourable circumstances this notion can be used to assign a braid
monodromy to a
family of polynomials.

\begin{description}

\item[Definition:]
A map $f:T\times\CC^2\to\CC$ is called a family of plane polynomials
admissible with
respect to a projection $p:\CC^2\to\CC$ if
\begin{enumerate}
\item
the restriction $f_t$ to each plane $\{t\}\times\CC^2$ is a polynomial,
\item
the zero divisor $Z_f=f\inv(0)$ and the singular values divisor $V_f$ are
branched covers
for the appropriate maps.
$$
\begin{array}{ccc@{\hspace*{1.5cm}}ccc}
Z_f & \to & T\times\CC^2 & V_f & \to & T\times\CC\\
&& \downarrow &&& \downarrow\\
&& T\times\CC &&& T
\end{array}
$$
\end{enumerate}
In this case the generalised braid monodromy of $V_f\subset T\times\CC$ is
called the
bifurcation braid monodromy of the family.

\end{description}

\begin{lemma}
\labell{a-two}
The bifurcation braid monodromy of the family $p_\l(x,y)=y^3-3\l y+2x$ is
generated by
the cube of the twist on the two singular values.
\end{lemma}

\proof
The divisor of singular values is given by the equation $\l^3=x^2$, hence the
bifurcation braid monodromy is the well known braid monodromy of a generically
projected simple cusp.
\qed

\begin{lemma}
\labell{a-one}
The bifurcation braid monodromy of the family $p_\l(x,y)=y^2-x^2+\l$ is the
full braid
group $\br_2$.
\end{lemma}

\proof
The divisor of singular values is given by the equation $x^2=\l$, hence the
bifurcation
braid monodromy is the well known braid monodromy of a vertical tangency
point on a
smooth double cover.
\qed

\begin{description}

\item[Definition:] The bifurcation braid monodromy group of a plane
polynomial $p_0$
with zero set a simple cover branched at $n$ points by a linear projection
$p:\CC^2\to\CC$ is the subgroup of $\br_n$ generated by the images of the
bifurcation
braid monodromy of all families of plane polynomials containing $p_0$ which are
admissible w.r.t.\ $p$.

\end{description}

\begin{prop}
\labell{deg-two}
The bifurcation braid monodromy of any generic polynomial deformation
equivalent to
$y^2-x^k$ is the full braid group $\br_k$.
\end{prop}

\proof
It suffices to consider the family $y^2-x^k+kx+\l$. Its singular value
divisor is given
by the equation $x^k-kx=\l$ of which the braid monodromy is as claimed.
\qed

Thus having dealt with the easiest cases we now want to investigate
polynomials with
branch degree three, which are in fact intimately related to the alternating
Artin systems considered in the first part of this paper.

\begin{lemma}
\labell{nat-basis}
The polynomial cover and its complement fibration for the polynomial
$y^3-3y+2x^k$ are
associated to the alternating Artin system of length $2k$ for a natural
choice of basis
of the fundamental group of the branch complement.
\end{lemma}

\proof
By straightforward computation the fibre at $x=0$ is $\CC^2$ with punctures at
$y=-\sqrt{3},0,\sqrt{3}$ which is regular, non-regular fibres occur at
$x^{2k}=1$ exactly
and along rays $x=r\zeta, r\in[0,1],\zeta$ primitive with $\zeta^{2k}=1$ the
points
$-\sqrt{3},0$ respectively $0,\sqrt{3}$ get closer and merge finally
according to
$\zeta^k=-1$ resp. $\zeta^k=1$.\\
So the elements of the star shaped basis are assigned alternating the twists
of the
intervals $[-\sqrt{3},0],[0,\sqrt{3}]$ which constitute the generator set
for an
alternating Artin system of length $2k$ generating $\br_3$.
\qed

Number the branch points of the polynomial cover given by the polynomial
$y^3-3y+2x^k$
according to increasing $arg$ starting with $x_1=1$.

Notice that the branch locus for the polynomial cover given by
$y^3-3p(x)y+2q(x)$ is
described by the equation $p^3(x)=q^2(x)$.\\

\begin{lemma}
\labell{S-deg}
The family $y^3-3y+2(x^k-\l),\l\in[-1,1]$ degenerates at $\l=\pm1$ only, branch
points are confined to straight rays, and even resp.\ odd indexed branch
points merge
at zero for $\l\to1$ resp.\ $\l\to -1$.
\end{lemma}

\proof
The branch points solve the equation
$$
\begin{array}{rccl}
(x^k-\l)^2=1 & \Leftrightarrow &
& x^k=\l+1\\
&&\vee & x^k=\l-1.
\end{array}
$$
The claim follows.
\qed

\begin{lemma}
\labell{S-mod}
The families $y^3-3y+2(x^k\pm1-\mu kx), \mu\in\CC$ small, have an
associated branch
locus divisor locally isomorphic to that of the family $y^2-x^k$.
\end{lemma}

\proof
The branch locus is given by the equation
$$
\begin{array}{cl}
&(x^k\pm1-\mu kx)^2=1\\
\Leftrightarrow & (x^k-\mu kx)(x^k-\mu kx\mp2)=0
\end{array}
$$
The corresponding divisor consists for small $\mu$ of a smooth unbranched
part and the
divisor associated to $y^2-x^k$.
\qed

\begin{lemma}
\labell{FF-deg}
The family $y^3-3(1-\l)y+2(x^k-i\l), \l\in[0,1]$ degenerates at $\l=1$
only, all branch
points are on a circle of modulus depending on $\l$, and pairs
$x_\nu,x_{\nu+1}$, where $\nu$ is even,
merge at $k$ distinct points for $\l\to 1$.
\end{lemma}

\proof
The branch points solve the equation
$$
\begin{array}{rccl}
(x^k-i\l)^2=(1-\l)^3 & \Leftrightarrow &
& x^k=i\l\pm\sqrt{(1-\l)^3}
\end{array}
$$
and one may check that $arg(x_\nu)$ is strictly increasing resp.\ decreasing
with $\l\to1$
for odd resp.\ even index.
\qed

\begin{lemma}
\labell{F-branch}
The family $y^3+\mu y+2(x^k-i-\mu),\mu\in\CC$ small, has an associated branch
locus divisor isomorphic to $k$ copies of the branch locus divisor of the
family $y^3-3\l y+2x$ locally at $\mu=0$.
\end{lemma}

\proof
The branch locus is given by the equation
$(x^k-i-\mu)^2=-\mu^3$. Hence up to invertible factors this equation reads at
each root $\a$ of the left hand side: $(x-\a)^2=\mu^3$. Since the number of
roots is $k$ the claim follows.
\qed

\begin{lemma}
\labell{F-deg}
There is a family of small deformations of the polynomial $y^3-2(x^k-i)$
such that
generically all branch points are simple except for a single double point.
\end{lemma}

\proof
At $\l=1$ the polynomial is of the form $y^3+2(x^k-i)$, hence the singular
values for
the projection are at $x^k=i$. On may choose an arbitrary one of these roots
say $\a$ and
define a perturbation $y^3-3\e(x-\a)y+2(x^k-i)$.
The singular values are now the zero locus of
$$
\e(x-\a)^3-(x^k-i)^2=:(x-\a)^2(\e(x-\a)-p_\a^2(x)).
$$
Assume there is another double root $x(\e)$ for the $\e$-family then
\begin{eqnarray}
\label{root}
\e(x(\e)-\a)-p_\a^2(x(\e)) & = & 0\\
\label{droot}
\wedge\hspace*{9mm} \e-2p_\a'(x(\e))p_\a(x(\e)) & = & 0
\end{eqnarray}
but on the other hand with equation (\ref{root}) also its derivative must
vanish:
\begin{eqnarray*}
0 & = &\frac\del{\del\e}\left(\e(x(\e)-\a)-p_\a^2(x(\e))\right)\\
& = & x(\e)-\a+\e x'(\e)-2p(x(\e))p'(x(\e))x'(\e)\\
& \stackrel{(\ref{droot})}{=} & x(\e)-\a+\e x'(\e)-\e x'(\e)\\
& = & x(\e)-\a.
\end{eqnarray*}
contrary to the assumption that $x(\e)$ is a root different from $\a$.
\qed

\begin{lemma}
\labell{F-mod}
There is a family of polynomials which contains in its interior a family
parameterised by $\l\in[0,1]$ such that $\l=0$ yields the polynomial
$y^3-3y+2x^k$ and for $\l\to1$ the family meets its only degeneration for which
all points branch points remain distinct except for the merging pair $x_1,x_2$.
Moreover the branch locus divisor has a cusp over $\l=1$.
\end{lemma}

\proof
We have to combine the families of the preceding lemmas into a family of two
complex parameters with $\a$ chosen to be the root at which $x_1,x_2$ merge.
Then the interval can be mapped to the parameter space in such a way as to
yield the desired properties.

The final claim follows from lemma \ref{F-branch}.
\qed

Our objective is now reached since the corresponding monodromies generate
the group $E_n$.

\begin{prop}
\labell{bif-mono}
The bifurcation braid monodromy group of the plane polynomials
$y^3-3y+2x^k$ contains a
subgroup in the conjugacy class of the standard isotropy group $E_n$.
\end{prop}

\proof
The case $k=1$ is lemma \ref{a-two}. For $k\geq2$ we consider the family
$y^3-3\ell(x)y+2(x^k-\mu)$
with $\ell(x)$ linear.
We compute the braid monodromy with respect to the natural choice of
geometric basis of
lemma \ref{nat-basis} with only a slightest move of the reference point
from the
origin to a point in
the sector defined by the rays of $x_{2k}$ and $x_1$.\\
It suffices to show that the generators $e_{12},e_{\nu,\nu+2}$ of $E_M$ are
contained in
the braid monodromy.
The triple twist $e_{12}$ is obtained by going around a degeneration as
given by the
family of lemma \ref{F-mod}.
The twists $e_{\nu,\nu+2}, \nu$ even resp.\ odd are generators for the full
braid group on
the even resp.\ odd indexed branch points. They are realised by an appropriate
deformation in the family of lemma \ref{S-mod}.
\qed

\section*{vanishing arcs and the Donaldson problem}

The objective of this last section is to shed some light on the Donaldson
problem of characterising the vanishing arcs among all isotopy classes of paths
in the base of a polynomial covering. We first sharpen the necessary criterion
of \cite{d}. Next we give a sufficient criterion in terms of the bifurcation
braid monodromy. Finally -- after a short digression to the situation for
finite coverings -- we make ends meet in favourable cases with the help of our
algebraic results.\\

First we give Donaldson's original definition of admissible paths, which we
call Coxeter admissible to distinguish them from our more restrictive
definition of Artin admissible paths:

\begin{description}

\item[Definition:]
A path in the base of a polynomial covering is called {\it Coxeter admissible},
if the monodromy of the finite cover along simple positive small loops around
the endpoints is one and the same transposition, if compared along the path.

\item[Definition:]
A path is called {\it Artin admissible}, if the braid monodromy assigned to
small simple positive loops around the endpoints is one and the same braid, if
compared along the path.

\item[Remark:] We could as well give a definition in terms of the geometry of
degeneration at the ends of a path. The concurrent pairs of points for the
degeneration to both ends of a Coxeter admissible path coincide. This pair
concurs even along isotopic paths in $\CC$ for the degeneration to both ends of
an Artin admissible path.

\end{description}
\pagebreak

\begin{lemma}
\labell{arc-adm}
A vanishing arc is Coxeter and Artin admissible.
\end{lemma}

\proof
A vanishing arc arises in a smoothing of a single ordinary double point.
So locally the divisor of critical values consists of only smooth
components without vertical tangents except for a single smooth component
locally isomorphic to a double cover branched at a single point.

Hence the vanishing arc is Artin admissible as claimed.
\qed

\begin{thm}
\labell{bif-orb}
The set of vanishing arcs for a plane polynomial $y^3-3y+2x^k$ is the orbit of a
single vanishing arc
under the action of the bifurcation braid monodromy group.
\end{thm}

\proof
Any vanishing arc is obtained by a degeneration of the reference polynomial
along an embedded path in the bifurcation complement with second endpoint a
generic point on the degeneration divisor.
Since versal families of the given polynomials are versal for the plane curve
singularity $y^3+x^k$, the degeneration
divisor is irreducible so we may assume within each class of such a path to have
chosen one which ends at a specified point. Hence a pair of vanishing arcs
defines an
element of the fundamental group by concatenation of the corresponding
paths. The braid associated to this loop maps one of the arcs into the other.
\qed

This result should be seen in contrast to the situation where instead of a
family of polynomial covers we consider (abstract) finite covers. Then the set
of corresponding vanishing arcs is much larger as it is even invariant under
the action of $S_{\cfami_{2k}}$
and coincides -- as remarked in \cite{d} -- with the set of Coxeter admissible
arcs.\\

We can finally characterise vanishing arcs as Artin admissible isotopy classes
in case the $x$-degree of our polynomial is sufficiently small:

\begin{thm}
The set of vanishing arcs for the polynomials $y^3-3y+2x^k$, $k=2,3$ is the
orbit of
the chord between the branch points $x_1,x_3$ under the bifurcation braid
monodromy and coincides with the set of Artin admissible paths.
\end{thm}

\proof
For the family of lemma \ref{S-mod} the chord between $x_1,x_3$ is a vanishing
arc and thm.\ \ref{bif-orb} then implies that its orbit under the bifurcation
braid monodromy group $E_{2k}$ is the set of vanishing arcs. By lemma
\ref{arc-adm} it is a subset of the Artin admissible arcs.

On the other hand if we are given an Artin admissible path then performing a
half twist on it does not change the monodromy of the polynomial covering,
hence it corresponds to a half twist in $S_{\afami_{2k}}$. This group coincides
with $E_{2k}$ in the given cases, so the half twist can actually be given as
$e_{13}^{\,e}$ with some $e\in E_{2k}$ by prop.~\ref{conclass}. We conclude
that the given path is the $e$-translate of our chord and hence is a vanishing
arc.
\qed

\begin{description}

\item[Remark:]
In general an Artin admissible isotopy class gives rise to a half twist
contained in $S_{\afami_{2k}}$ and it is an open question whether
$S_{\afami_{2k}}$ equals $E_{2k}$ and whether its half twists are contained in
a single $E_{2k}$ conjugation class.

\end{description}

\pagebreak


\begin{thebibliography}{CBW}

\bibitem[Ar]{ar} E. Artin: {\sl Theory of braids},
Ann. Maths. 48 (1947), 101--126

\bibitem[AK]{ak} D. Auroux, L. Katzarkov:
{\sl Branched coverings of ${\bf C}{\rm P}\sp 2$ and invariants of symplectic
   4-manifolds},
Invent. Math. 142
(2000), 631--673

\bibitem[BW]{bw} J. Birman, B. Wajnryb:
{\sl $3$-fold branched coverings and the mapping class group of a surface},
in Geometry and topology (College Park, Md., 1983/84), ed. J. Alexander, J.
Harer, 24--46,  Lecture Notes in Math., 1167,
Springer, Berlin-New York, 1985.

\bibitem[CW]{cw}
F. Catanese, B. Wajnryb:
{\sl The fundamental group of generic polynomials},
Topology 30 (1991), no. 4, 641--651

\bibitem[CS]{cs} D.~Cohen, A.~Suciu:
{\sl The braid monodromy of plane algebraic curves and hyperplane arrangements},
Comment.\ Math.\ Helv.\ 72 (1997), 285--315

\bibitem[CP]{cp}
J.\ Crisp, L.\ Paris:
{\sl The solution to a conjecture on the subgroup generated by the squares of
the generators of an Artin group}, Invent.\ 145 (2001), no.\ 1, 19--36

\bibitem[D\"o]{do} A. D"orner:
{\sl Isotropieuntergruppen der artinschen Zopfgruppen},
Bonner Mathematische Schriften 255.  Universit"at Bonn, Mathematisches
Institut, Bonn, 1993.

\bibitem[Do]{d}
S. K. Donaldson:
{\sl Polynomials, vanishing cycles and Floer homology},
in Mathematics: frontiers and perspectives, 55--64,
AMS., Providence, 2000

\bibitem[Ha]{ha} V.~L.~Hansen:
{\sl Braids and Coverings},
Cambridge Univ.\ Press (1989),
London Math.\ Soc.\ Student Texts {\bf{18}}

\bibitem[Hu]{hu} A. Hurwitz: {\sl \"Uber Riemann'sche Fl\"achen mit gegebenen
Verzweigungspunkten},
Math. Ann. 39 (1891), 1--61

\bibitem[Kl]{klui} P. Kluitmann: {\sl Isotropy Subgroups of Artin's Braid
Group}, preprint 1991

\bibitem[KT]{kt} V. S. Kulikov, M. Teicher:
{\sl Braid monodromy factorisations and
diffeomorphism types}, izv. Math. 64, 311--341 (2000)

\bibitem[Lo]{lo}
E.~Looijenga:
{\sl The complement of the bifurcation variety of a simple singularity},
Invent. Math. 23 (1974), 105--116.

\bibitem[Mo]{mo}  B. Moishezon:
{\sl Stable branch curves and braid monodromies}
in Algebraic geometry (Chicago, Ill., 1980),
Lecture Notes in Math. 862, 107--192, Springer, Berlin, 1981

\bibitem[MT]{mt} B. Moishezon, M. Teicher:
{\sl Braid group technique in complex geometry. I. Line arrangements in ${\bf
   C}{\rm P}\sp 2$},
in Braids (Santa Cruz, CA, 1986), Contemp. Math. 78, 425--555, AMS,
Providence, 1988

\end{thebibliography}
\end{document}